\documentclass[10pt,a4paper]{article}
\usepackage[latin1]{inputenc}
\usepackage{amsmath}
\usepackage{amsfonts}
\usepackage{amssymb}
\title {A Proof of the Quadratic Reciprocity Law}
\author{Math Dicker,  \small {Open University of the Netherlands}}

\begin{document}

\maketitle
\begin{abstract}
A proof of the \textit{Quadratic reciprocity Law} is presented using a \textit{Lemma of Gauss}, \textit{the theory of finite fields} and the \textit{Frobenius automorfism}.
\end{abstract}
\section{Introduction.}
Let $p$, $q$ be distinct odd prime numbers and let $e$ denote the order of $q$ in  $ \mathbb{F}_{p}^{*} $. The Frobenius automorphism $x\rightarrow x^q$  in the  field $ \mathbb{F}_{q^{e}} $ is here denoted by $ \varphi_{q} $. Because $p$ divides $q^e-1$, the cyclic group $ \mathbb{F}_{q^{e}}^{*} $ contains a primitive $p$-th root of unity to which we refer by $\theta$. If we specify  $f(x)=1+x+x^2....+x^{p-1}$  then $f(\theta^{k})=0$ for k with gcd(k,p)=1, otherwise $f(\theta^{k})=p$. We denote by $\delta(x_{1},x_{2},x_{3},.... , x_{p})$  the  determinant of the $p$-square matrix with the entry in the $i$th row and $j$th column equal to $ (x_{j})^{i} $.
\section{The Quadratic Reciprocity Law.}
In particular, $\delta(1,\vartheta,\vartheta^2,.... ,\vartheta^{p-1})$  is the following determinant: 
\[ \\ \delta(1,\vartheta,\vartheta^2,.... , \vartheta^{p-1})=       
\left|
\begin{array}{ccccccccc} 
1&\vartheta&\vartheta^2&\vartheta^3&\cdots &\vartheta^{p-1}\\
1&\vartheta^2&\vartheta^4&\vartheta^6&\cdots &\vartheta^{2(p-1)}\\
1&\vartheta^3&\vartheta^6&\vartheta^9&\cdots &\vartheta^{3(p-1)}\\
\vdots&\vdots &\vdots &\vdots &\ddots &\vdots \\
1&\vartheta^i&\vartheta^{2i}&\vartheta^{3i}&\cdots &\vartheta^{i(p-1)}\\
\vdots&\vdots &\vdots &\vdots &\ddots &\vdots \\
1&1&1&1&\cdots &1\\
\end{array}
\right|
\]
\\
Let M be the matrix corresponding with this determinant.
\\[2ex]
\textbf{Theorem}
\\[2ex]
1] $\delta(1,\vartheta,\vartheta^2,.... , \vartheta^{p-1})^2=p^*p^{p-1}$ with $p^*=(-1)^{\frac{p-1}{2}}p$
\\[2ex]
2] $\varphi_{q} (\delta(1,\vartheta,\vartheta^2,.... , \vartheta^{p-1}))  = (\dfrac{q}{p})  \delta(1,\vartheta,\vartheta^2,.... , \vartheta^{p-1})$
\\[2ex]
From 1], 2] and using $\varphi_{q} (x)=x $ $\Longleftrightarrow$ $x\in$ $ \mathbb{F}_{q} $ and Euler's criterion, it follows:
\\[2ex]
$ (\dfrac{p*}{q})=1 \Leftrightarrow (\dfrac{q}{p})=1 $ or $ (\dfrac{p}{q})(\dfrac{q}{p})=(-1)^{\frac{p-1}{2}\frac{q-1}{2}} $ 
\\[2ex]
\textbf{Proof}
\\[2ex]
ad 1] Consider the  matrixproduct: 
\[ \\ M^{T}M=       
\left|
\begin{array}{ccccccccc} 
p&0&0&0&\cdots &0\\
0&0&0&0&\cdots &p\\
0&0&0&\cdots &p&0\\
\vdots&\vdots &\vdots &\vdots &\ddots &\vdots \\
0&0&0&p&\cdots &0\\
\vdots&\vdots &\vdots &\vdots &\ddots &\vdots \\
0&p&0&0&\cdots &0\\
\end{array}
\right|
\]
because $ M^{T}M $ = $\left( f(\vartheta^{(i+j-2)})\right) $ for row $i=1,2,.,p$ and column $j=1,2,., p$.
\\[2ex]
ad 2] Consider the residue classes of $ \mathbb{F}_{p}^{*} $  represented by the following half systems:
H $=1,2,3,...,\frac{(p-1)}{2}$ and -H $=-1,-2,-3,...,-\frac{(p-1)}{2}$.\
We introduce the  function $\varrho_{q}$ which is connected as we shall see in a moment, to the Frobenius automorfism;   The function $\varrho_{q} :  \mathbb{F}_{p}^{*}  \longrightarrow  \mathbb{F}_{p}^{*} $ is defined by $x \longrightarrow qx$;  the result of $\varrho_{q}$ is a permutation of $ \mathbb{F}_{p}^{*} $. If we denote by $\mu $
the number of elements in the set S, with $S=\{x|x\in H ,\varrho_{q}(x)\notin H \} $, then  $(\dfrac{q}{p})=(-1)^{\mu}$
(lemma of Gauss)\cite{gauss}.
\\[2ex]
Important for us is that the permutation $\varrho_{q}$ working on $ \mathbb{F}_{p}^{*} $  is  the result of  $\mu$ interchanges, leaving aside a multiple of $2$. This can be grasped as follows: define the permutation $\pi$ on $ \mathbb{F}_{p}^{*} $ :\\
$\pi=\left[ \Pi_{y\in S}(\varrho_{q}(y),-\varrho_{q}(y))\right]  \varrho_{q}$; by (i,j) is denoted the permutation which interchanges i and j. The permutation $\varrho_{q}$ and $\pi$ originate from each other by $\mu$ interchanges.
\\[2ex]
The permutation $\pi$ has the following properties: $\pi(x)=-\pi(-x) $ ; $ \pi(H)=H$ and $\pi(-H)=-H$. Hence the identity permutation originates from $\pi$ by an even number of paired interchanges on H and -H.
\\[2ex]
We have: $\delta(1,\vartheta,\vartheta^2,.... , \vartheta^{p-1})$ =
$\delta(1,\vartheta,\vartheta^2,.., \vartheta^{\frac{p-1}{2}},\vartheta^{-\frac{(p-1)}{2}},..,\vartheta^{-2},\vartheta^{-1})$.
\\[2ex]
The application of the Frobenius automorfism and the consideration of the above mentioned properties  of $\varrho_{q}$, results in:
\\[4ex]
$\varphi_{q}$ $(\delta(1,\vartheta,\vartheta^2,.., \vartheta^{\frac{p-1}{2}},\vartheta^{-\frac{(p-1)}{2}},..,\vartheta^{-2},\vartheta^{-1}))$ =
\\[2ex]
$\delta(1,\vartheta^{q},\vartheta^{2q},.., \vartheta^{\frac{p-1}{2}q},\vartheta^{-\frac{(p-1)}{2}q},..,\vartheta^{-2q},\vartheta^{-1q})$=
\\[2ex]
$ (-1)^{\mu} $  $\delta(1,\vartheta,\vartheta^2,.., \vartheta^{\frac{p-1}{2}},\vartheta^{-\frac{(p-1)}{2}},..,\vartheta^{-2},\vartheta^{-1})$=
\\[2ex]
$ (\dfrac{q}{p})$  $\delta(1,\vartheta,\vartheta^2,.., \vartheta^{\frac{p-1}{2}},\vartheta^{-\frac{(p-1)}{2}},..,\vartheta^{-2},\vartheta^{-1})$
\\[2ex]
\section{An example of $\varrho_{q}$.}
The crucial point of this proof is that the operation of $\varphi_{q}$ on the determinant $\delta(1,\vartheta,\vartheta^2,.... , \vartheta^{p-1})$ results in $ \mu $ interchanges of the columns, leaving aside a multiple of $2$ interchanges. This is caused by the properties of the connected function $\varrho_{q}$.\\[2ex]
To illustrate, we consider the case where $p=13$ and $q=5$. The residue classes of $ \mathbb{F}_{13}^{*} $ are represented by:
$1,2,3,4,5,6,-6,-5,-4,-3,-2,-1$. H contains the first six classes. The function $\varrho_{5}$ results in the following permutation: $5,-3,2,-6,-1,4,-4,1,6,-2,3,-5$. $\mu $ is 3; hence  $(\dfrac{5}{13})=-1$ [lemma of Gauss]; $S=\left\lbrace {2,4,5}\right\rbrace$. The permutation $\pi$ is: $5,3,2,6,1,4,-4,-1,-6,-2,-3,-5$. The permutation $\pi$ can be transformed into the identity permutation by the following six paired interchanges (-2,-3) (2,3) (-4,-6) (4,6) (-5,-1) (5,1)

Math Dicker, Kerselarenstraat 36\\
3700 Tongeren, Belgium\\
louis.dicker@telenet.be; math.dicker@ou.nl

\end{document}